\documentclass[12pt,notitlepage,twoside,a4paper]{amsart}
 \usepackage{amsfonts} 
\usepackage[latin1]{inputenc}
\usepackage[cyr]{aeguill}
\usepackage{xspace}
\usepackage[polutonikogreek,english]{babel}

\usepackage{amsmath,amssymb,enumerate}

\usepackage{epsfig,fancyhdr,color}

\usepackage{amssymb}
\usepackage{amsmath,amsthm}
\usepackage{latexsym}
\usepackage{amscd}
\usepackage{psfrag}
\usepackage{graphicx}
\usepackage[latin1]{inputenc}
\usepackage[all]{xy}
\usepackage[mathcal]{eucal}

\definecolor{NoteColor}{rgb}{1,0,0}

\renewcommand{\textsc}{\textcolor{red}}

\newtheorem*{theorem 1}{\rm\bf Proposition 1}
\newtheorem*{theorem 2}{\rm\bf Proposition 2}

\theoremstyle{definition}

\theoremstyle{remark}

\def\interieur#1{\mathord{\mathop{\kern 0pt #1}\limits^\circ}}

\title[On Delisle's geographical projection]{On Delisle's geographical projection}

\author{Charalampos Charitos and Athanase Papadopoulos}
\address{Charalampos Charitos, Agricultural University Athens, Mathematics Laboratory, Iera Odos 75, 11855 Athens, Greece, email: bakis@aua.gr;  Athanase Papadopoulos,  Universit{\'e} de Strasbourg and CNRS,
7 rue Ren\'e Descartes,
 67084 Strasbourg Cedex, France, email: papadop@math.unistra.fr}
\makeindex
 
\date{\today}

\begin{document}

  \begin{abstract}  Joseph-Nicolas Delisle was one of the most important scientists at the Saint Petersburg Academy of Sciences during the first period when Euler was working there. Euler was helping him in his work on astronomy and in geography. In this paper,
  Delisle's geographical projection is presented and Euler's study of this projection is
explained, highlighting some important mathematical points, in particular on the metric geometry of surfaces.

The final version of this paper will appear in the book \emph{Mathematical Geography in the Eighteenth Century: Euler, Lagrange and Lambert}, ed. 
Renzo Caddeo and Athanase Papadopoulos, Springer International Publishing, 2022.
      
       
      \medskip
      

     \end{abstract}
          
  \maketitle
  
  \tableofcontents

\section{ Introduction}

\begin{figure}[htbp]

\centering
\includegraphics[width=12.5cm]{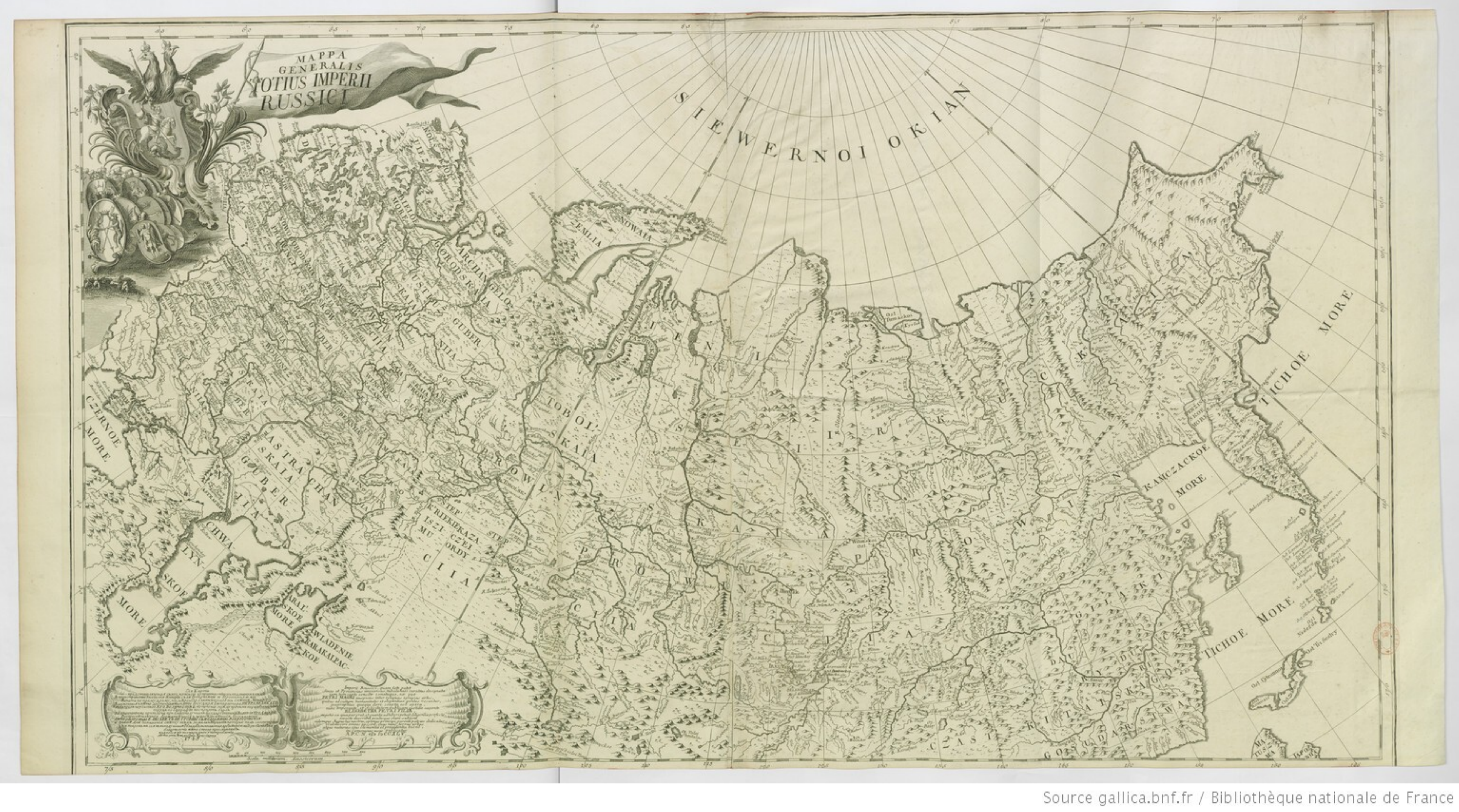}
\caption{\small{A map of the Russian Empire, drawn using Delisle's method:\index{map!Russian Empire}\index{general map of the Russian Empire} this is the last map of Euler's  \emph{Atlas Russicus}, published in Saint Petersburg in 1745 \cite{Atlas-Russicus}. Bibliothèque nationale de France, Département Cartes et Plans}} \label{Fig:Empire}
\end{figure}

The famous French geographer Joseph-Nicolas Delisle,\index{Delisle, Joseph-Nicolas}  who was Leonhard Euler's colleague and collaborator at the Saint Petersburg Academy of Sciences, introduced a projection from the sphere onto the Euclidean plane which became known as Delisle's geographical projection. In Figure \ref{Fig:Empire}, we have reproduced a map of the Russian Empire drawn under the direction of Euler{The project, before Euler took it over, was directed by Delisle}, that uses Delisle's projection.\index{map!Delisle}\index{projection!Delisle}\index{Delisle map}\index{Delisle projection} 
This projection shares several properties of the conical projection which was used in Greek Antiquity, although the two projections are different. The conical projection is 
 described in  Chapters 21 and 24 of Book I of Ptolemy's\index{Ptolemy, Claudius} \emph{Geography}  \cite{Ptolemy-geo}. It is obtained by first projecting the surface of the Earth onto a cone tangent to it along a certain parallel, in such a way that each meridian is sent to the line in which the plane containing it intersects the cone. This cone is then unfolded onto a plane. In this way, parallels are sent to concentric circles and meridians to straight lines with a common intersection point (which is not the North pole). Distances are preserved on the parallel that we started with.

  Ptolemy\index{Ptolemy, Claudius}, who was aware of the fact that it is not possible for a geographical map to preserve proportions of all distances, used a conical projection in which these proportions are preserved along two special parallels, namely, the parallel passing through the island of Thule (the farthest northern location mentioned in the \emph{Geography}, and in some other geographical works of Greek antiquity)\footnote{Some authors conjecture that Thule is an island in Norway, some others that it is Groenland, and there are other possibilities} and the equator.  He then discussed the  corrections that have to be made in the region between these two parallels in such a way that the distortion there is optimal. In Delisle's projection, which is the subject of the present chapter, like in Ptolemy's projection, proportions of distances are also preserved along two chosen parallels. In particular, in Delisle's map of the Russian Empire,  these parallels are those that bound this Empire. Like in Ptolemy's conical projection, Delisle's map is constructed, in each case, so that the distortion is optimized between these parallels.

We note incidentally that the preservation of ratios of distances on some special parallel is also a property of the cylindrical projection of Marinus of Tyre,\index{Marinus of Tyre} see the discussion by Neugebauer\index{Neugebauer, Otto} in \cite[p. 1037--1039]{Neu1948} where this author
mentions a cylindrical projection\index{cylindrical projection}\index{map!cylindrical}  by Marinus in which ratios of distances are preserved along all the meridians and along the parallel passing through the island of Rhodes.

\begin{figure}
\begin{center}
\includegraphics[scale=0.5]{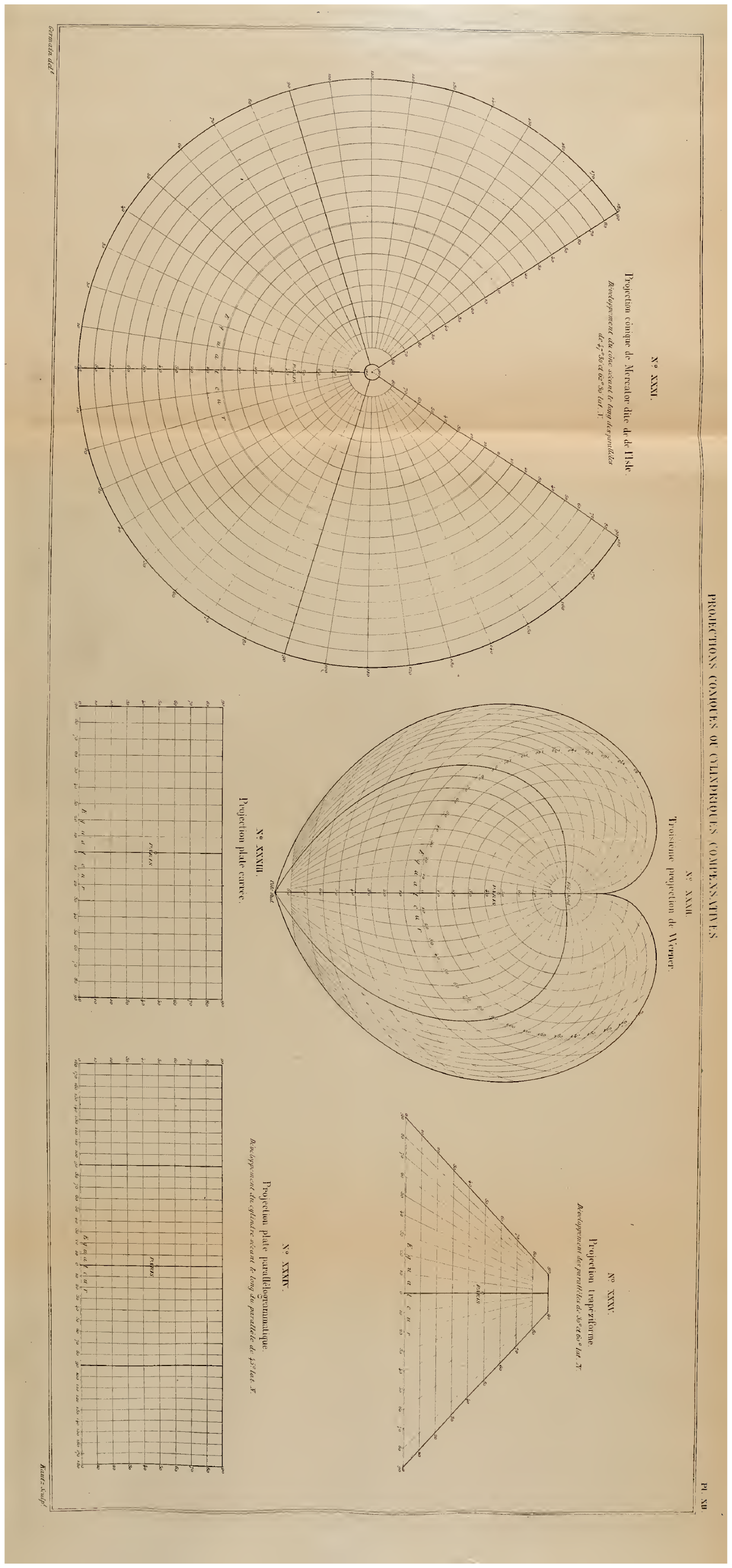} 
\end{center}
\label{fig:MD}
\caption{{\small The so-called Mercator--Delisle projection, drawing from  \cite{Germain}}}
\end{figure}

We mention now another projection where one starts with a cone which is not tangent  to the sphere but which intersects it in two parallels. This projection is more closely related to the one of Delisle. It was used by Gerardus Mercator\index{Mercator, Gerardus (Gheert Cremer)} in his 1544 construction of the map of Europe; see \cite[p. 178-179]{Germain} where the author considers that Delisle's map is in fact the  same as the one of Mercator. In any case, there is no doubt that Delisle was familiar with Ptolemy's maps, with other maps of Antiquity, and with Mercator's maps.
  In both projections (Ptolemy and Mercator), parallels are sent to concentric circles and meridians to straight lines intersecting at a point which is not the image of the North pole. Delisle's projection\index{map!Delisle}\index{projection!Delisle}\index{Delisle map}\index{Delisle projection} satisfies the same properties. In Figure \ref{fig:MD}, we have reproduced a drawing from \cite{Germain} that represents what the author calls the Mercator--Delisle projection.

Finally, we note that in Lambert's memoir on geography that is translated in the present volume, the latter describes also a conical projection in which two circles of latitude  (that is, two parallels)  are mapped such that on the images of these parallels,  the proportions are preserved.

Euler, in his memoir titled\index{map!Delisle}\index{projection!Delisle}\index{Delisle map}\index{Delisle projection} \emph{De proiectione geographica De Lisliana in mappa generali imperii russici usitata} (On Delisle's geographical projection used for a general map of the Russian Empire), published in 1778  (see \cite{Euler}, translated from the Latin in the present volume), gave a mathematical description of Delisle's projection.
Our goal in this chapter is to present, using modern mathematical notation, Euler's description of this projection, highlighting certain interesting mathematical points.

 Delisle's\index{map!Delisle}\index{projection!Delisle}\index{Delisle map}\index{Delisle projection} projection satisfies the following requirements:

\begin{itemize}
\item All meridians are represented by straight lines.

\item The projection preserves the degrees of latitudes, i.e. the projection is faithful along the meridians.

\item Meridians and parallels intersect perpendicularly.
\end{itemize}

All straight lines which are images of meridians intersect at a common point (we shall linger on this below). Such a map cannot
preserve at every point the ratio of the length corresponding to a single degree
of the parallel passing through this point to the length corresponding to a single degree of the meridian.
 Euler thus writes (see
(\cite{Euler} \S 5) that the
following general question is of great importance: \textquotedblleft \textit{%
In what way should the meridians be arranged with respect to the parallels so that for the whole extent of the map, 
the deviation from the ratio that the degrees of longitude and latitude which [the meridians and parallels] have among themselves on the sphere is the smallest possible?}%
\textquotedblright

In Delisle's\index{map!Delisle}\index{projection!Delisle}\index{Delisle map}\index{Delisle projection} geographical projection, two special parallels are chosen, 
along which  the projection is an equidistant map, that is, the
proportion between the degrees of longitude and latitude is
preserved. Furthermore, Delisle\index{map!Delisle}\index{projection!Delisle}\index{Delisle map}\index{Delisle projection} discovered that if these parallels are equidistant from the central
parallel and one from the southernmost parallel and the other from the
northernmost parallel, then the deviation of the
geographical map is nowhere significant. Thus, another question is \textit{to
detect these two parallels mentioned above so that even the largest errors
that may arise in the map are the smallest possible.}

One may also find mathematical comments on Delisle's projection, also based on 
Euler's work, in \cite{Steffens}, \S 1.1 titled \emph{Euler's analysis on Delisle's map} and in \cite{Germain}, Chapter VI, \S 5.

\section{Construction of Delisle's projection}

In the rest of this chapter,\index{map!Delisle}\index{projection!Delisle}\index{Delisle map}\index{Delisle projection} following Euler \cite{Euler}, the figure of the Earth is considered to be 
spherical.  The fact that the Earth was known to be rather a spheroid is not taken into account, since the difference between this spheroid and a sphere will not be visible on a geographical map. Thus, all meridians are great
circles, or semi-circles, depending n the context, of the same length.

 In what follows, the sphere representing the Earth is denoted by $S$ and $\delta$ will denote the length of the radius of every meridian.

For points $X,$ $Y$ on $S$ an arc in $S$ joining $X$ and $Y$ will be denoted
by $XY.$ In practice, this arc will be either an arc
of a meridian or an arc of a parallel
of $S$ and in each case we shall determine it precisely. On the other hand, the segment of straight line joining $X$ and $Y$ in the ambient Euclidean
space $E^{3}$ will be denoted by $\overline{XY}$ and will be referred to as a 
\textit{line segment}. Also, if $X$ and $Y$ are points on a meridian of $S$, then by an abuse of notation
we shall denote by $XY$ the meridian (that is, the great circle or semi-circle) determined by $X$ and $Y.$
Finally for points $X,$ $Y$ on $S$, the distance of these points in $S$ will
be referred to as the \textit{Earth distance} between them while the distance of these points
in $E^{3}$ will be referred to as the \emph{Euclidean distance} between them.

\begin{figure}
\begin{center}
\includegraphics[scale=1]{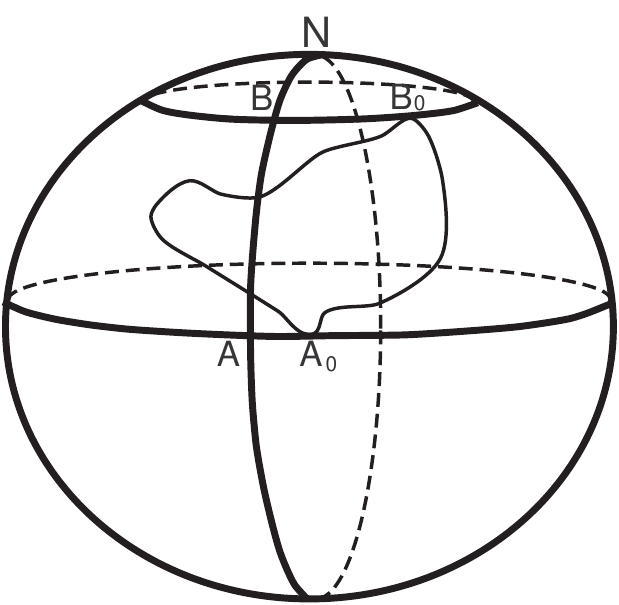} 
\end{center}
\caption{\small The southernmost and the northernmost parallels that bound the Russian
Empire and the principal meridian $AB.$}
\label{f:1}
\end{figure}

Let $P$ be a point of latitude $p$ on the sphere $S$ of radius $\delta$ representing the Earth. The length of one degree of meridian is equal to $\delta$. 
 The length of one degree of longitude on the
parallel passing through $P$ is equal to $\delta \cos p.$ Thus, the ratio of one degree of longitude on the parallel passing trough $P$ to the length of one degree of latitude  remains
constant on each meridian and equal to $\cos p.$

Let $A_{0}$ and $B_{0}$ be respectively the southernmost and the northernmost points of
the Russian Empire. The parallels passing through these points will be referred to as the \textit{southernmost}
and the \textit{northernmost borders} respectively, see Figure \ref{f:1}. The
latitude of $A_{0}$ is $a=40^{\mathrm{o}}$ and the latitude of $B_{0}$ is $b=70^{\mathrm{o}}$
approximately. We consider now an arbitrary meridian intersecting the
southernmost and the northernmost parallels at the points $A$ and $B$
respectively, see Figure  \ref{f:1}.

On this meridian $AB,$ which will be referred to as the \textit{principal
meridian},\index{principal meridian}\index{meridian!principal} we consider the points $P$ and $Q$ whose latitudes are $p=50^{\mathrm{o}}$ and 
$q=60^{\mathrm{o}}$ respectively. Finally on the parallels passing through the points 
$P$ and $Q$ we consider points $P_{1}$ and $Q_{1}$ respectively such that
the length of the arcs $PP_{1}$ and $QQ_{1}$ (contained in the parallels) are $\delta \cos p$ and $\delta \cos q$ respectively. In other words, the
longitude of both points $P_{1}$ and $Q_{1}$ (measured from the principal
meridian) is one degree, see Figure   \ref{f:2}.

Since the length of a degree of meridian is small for our purposes, the lengths of the line segments $\overline{PP_{1}}$ and $\overline{QQ_{1}}$ are almost equal to those of the small arcs $PP_{1}$ and $QQ_{1}$. That is, we
may consider that the lengths of the Euclidean segments $\overline{PP_{1}}$ and $\overline{QQ_{1}}$  are $%
\delta \cos p$ and $\delta \cos q$ respectively. These segments may be considered perpendicular to the meridian $AB$.

We wish to define a
projection $f:U\rightarrow E^{2}$ satisfying the requirements of the
introduction, where $U$ is an open subset of $S$ containing the Russian
Empire and $E^{2}$ is the Euclidean plane. For this, we shall first consider four
specific points in $E^{2}$ which will serve as images of $P,$ $P_{1},$ $Q,$ $%
Q_{1}$ by $f.$ These points will be denoted by $P^{\prime },$ $P_{1}^{\prime
},$ $Q^{\prime },$ $Q_{1}^{\prime }.$ Gradually, from these points a network
consisting of the images of meridians and parallels will be constructed in $%
E^{2}$ and thus the projection is defined.

We first consider two points $P^{\prime },$ $Q^{\prime }$ in $E^{2}$ with
distance $|P^{\prime }Q^{\prime }|$ $=\delta (q-p);$ that is, the Euclidean
distance $|P^{\prime }Q^{\prime }|$ in the plane $E^{2}$ is equal to the distance
of the points $P$ and $Q$ in $S$. The segment $P^{\prime }Q^{\prime
}\subset E^{2}$ will be the image of the arc of meridian $PQ\subset S.$ Now,
we consider points $P_{1}^{\prime }$ and $Q_{1}^{\prime }$ in $E^{2}$ such
that:

\begin{figure}
\begin{center}
\includegraphics[scale=0.9]{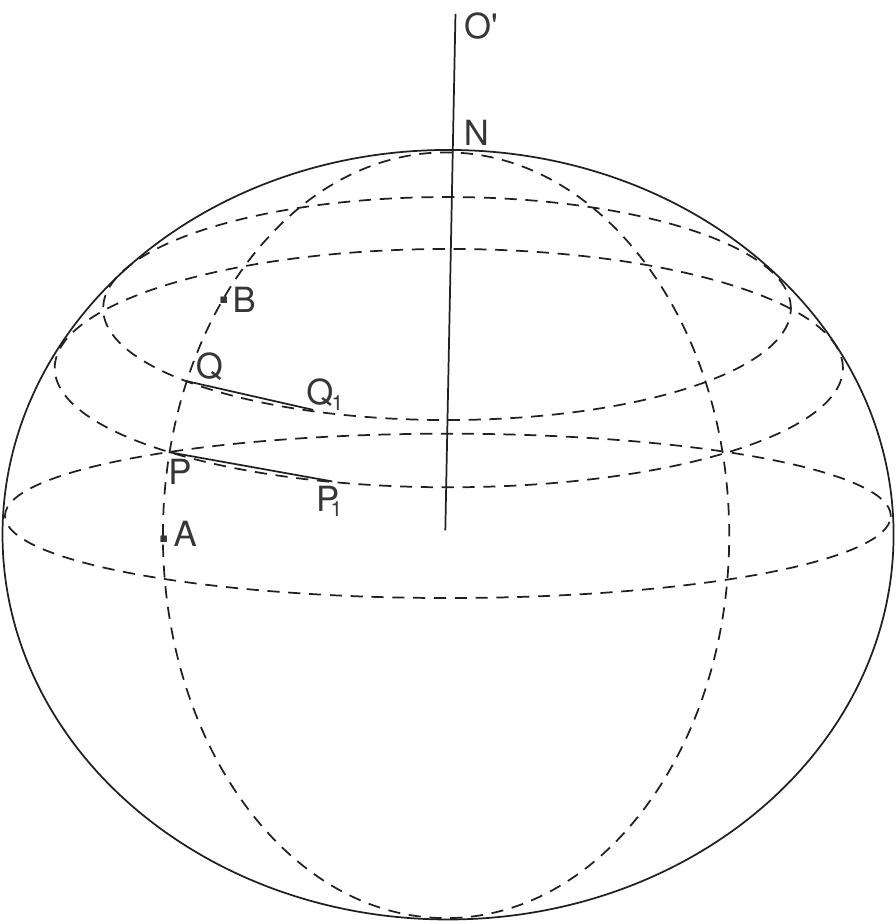} 
\end{center}
\caption{\small The principal meridian and the line segments $PP_{1}$ and $QQ_{1}.$}
\label{f:2}
\end{figure}

(1) $P^{\prime }P_{1}^{\prime }$ and $Q^{\prime }Q_{1}^{\prime }$ are
segments in $E^{2}$ perpendicular to $P^{\prime }Q^{\prime };$

(2) $|P^{\prime }P_{1}^{\prime }|$ $=\delta \cos p$ and $|Q^{\prime
}Q_{1}^{\prime }|$ $=\delta \cos q;$ that is, the distance $|P^{\prime
}P_{1}^{\prime }|$ is the length of the arc $PP_{1}$ contained in the
parallel passing through $P.$ Similarly for the distance $|Q^{\prime
}Q_{1}^{\prime }|$; see Figure \ref{f:2}. By definition, the points $P^{\prime },$ $%
P_{1}^{\prime },$ $Q^{\prime },$ $Q_{1}^{\prime }$ will be the images of $P,$
$P_{1},$ $Q,$ $Q_{1}$ respectively and the segments $P^{\prime }Q^{\prime }$
and $P_{1}^{\prime }Q_{1}^{\prime }$ in $E^{2}$ will be the  images of the arcs of meridian $PQ$ and $P_{1}Q_{1}.$

We let $O^{\prime }$ be the intersection point of the lines containing the segments $P^{\prime }Q^{\prime }$ and $%
P_{1}^{\prime }Q_{1}^{\prime }$ in $E^{2}$. We shall
 compute the length of $O^{\prime }P^{\prime }.$ From Figure 4, we
have:

\begin{equation*}
\frac{|P^{\prime }P_{1}^{\prime }|-|Q^{\prime }Q_{1}^{\prime }|}{|P^{\prime
}Q^{\prime }|}=\frac{|P^{\prime }P_{1}^{\prime }|}{|P^{\prime }O^{\prime }|}
\end{equation*}%
or, 
\begin{equation*}
\frac{\delta (\cos p-\cos q)}{\delta (q-p)}=\frac{\delta \cos p}{|P^{\prime
}O^{\prime }|}.
\end{equation*}%
Therefore,%
\begin{equation}
|P^{\prime }O^{\prime }|\text{ }=\frac{\delta (q-p)\cos p}{\cos p-\cos q}. 
\tag{1}  \label{1}
\end{equation}%
This formula in degrees takes the form 
\begin{equation}
|P^{\prime }O^{\prime }|\text{ }=\frac{(q-p)\cos p}{\cos p-\cos q}.  \tag{2}
\label{2}
\end{equation}

Note that in this discussion, the fact that the distance between two points $Z,$ $%
W\in E^{2}$ is expressed in degrees $\theta$ means that $|ZW|$ is equal to $%
\delta \theta.$ In other words, $|ZW|$ is the length of an
 arc of meridian of $S$ of $\theta $ degrees.

Setting $p=50^{\mathrm{o}},$ $q=60^{\mathrm{o}}$ in (\ref{2}) we find that $|P^{\prime
}O^{\prime }|$ is the length of an arc of $45^{\mathrm{o}}1^{\prime }$ in a
meridian of $S.$ Therefore, since the point $P$ in $S$ is at distance $50^{\mathrm{o}}$
from the equator with respect the Earth distance of $S,$ we deduce that if $%
O^{\prime }$ is placed in the ambient space $E^{3}$ of $S$, namely, on the
vertical straight line passing trough $N,$ then $O^{\prime }$ lies
beyond the Earth's North pole $N$ at a Euclidean distance $\delta \cdot
45^{\mathrm{o}}1^{\prime }$ from $P$ or at a Euclidean distance $\delta \cdot
95^{\mathrm{o}}1^{\prime }$ from the equator.

\begin{figure}
\begin{center}
\includegraphics[scale=0.7]{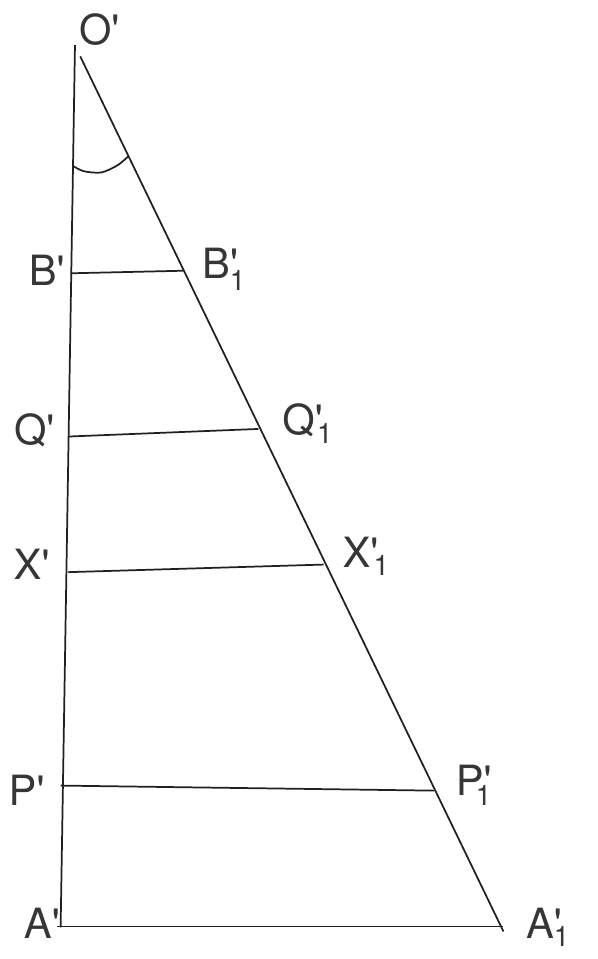} 
\end{center}
\caption{\small The image of the principal meridian in the Euclidean plane $E^{2}.$}
\label{f:4}
\end{figure}

Now we can construct Delisle's projection\index{map!Delisle}\index{projection!Delisle}\index{Delisle map}\index{Delisle projection} in $E^{2}.$ The images of
all the meridians will be straight lines emanating from $O^{\prime }.$
Considering the line containing $O^{\prime }P^{\prime }$, we draw a circle $%
S_{P^{\prime }}$ of center $O^{\prime }$ and radius $O^{\prime }P^{\prime }$
which is divided into parts equal to $\delta \cos p,$ which is the length of
one degree of the parallel passing through $P.$ The lines led from $%
O^{\prime }$ and passing through each of the points in the subdivision will
give all the meridians that are drawn on the map. This being done, 
 together with the circle S$S_{P^{\prime }}$, we draw  all the circles of center 
$%
O^{\prime }$ which are $\delta $  distant  apart. These circles are the
images of parallels of $S$ and the distance between two consecutive circles is one degree of
latitude (Figure \ref{f:4}). In this way a projection $f$\ satisfying all the requirements of
the introduction is defined.

  In the next section
calculating the distortion of $f$ we will deduce that the
 choice of $P$ ($p$ = 50$^\mathrm{o}$) and $Q$ ($p$ = 60$^\mathrm{o}$)
choice of $P$ and $%
Q $ is extremely close to the best choice of parallels along which the ratio
of degrees of longitude and latitude is correct.

\section{The distortion of Delisle's projection}

In this section,\index{map!Delisle}\index{projection!Delisle}\index{Delisle map}\index{Delisle projection} in a more general context and following Euler's mathematical
analysis, we will estimate the distortion from reality of Delisle's
projection in order to be able in the next section to detect the specific
points $P$ and $Q$ on $S$ mentioned above and apply the whole study to the
case of the map of the Russian Empire.

First we will compute the angle $\omega =P^{\prime }O^{\prime }P_{1}^{\prime
}$ corresponding to one degree of longitude in the map, see Figure 3. We may assume that the length of the circular arc of radius $\vert O'P'\vert$ and of angle $\omega$ is approximately equal to  
$|P^{\prime }P_{1}^{\prime }|.$ Therefore, from (\ref{1}) the angle $\omega $ in radians 
is equal to

\begin{equation}
\omega =\frac{|P^{\prime }P_{1}^{\prime }|}{|P^{\prime }O^{\prime }|}=\frac{%
\cos p-\cos q}{\alpha (q-p)}  \tag{3}  \label{3}
\end{equation}%
where the factor $\alpha =\frac{\pi }{180}=0.01745329$ measures one degree
in radians. With this formula we may see that $\omega =49^{\prime }6^{\prime
\prime }$ in degrees.

Now given that the latitude of $P$ is equal to $p,$ we consider the point $%
N^{\prime }$ in $O^{\prime }P^{\prime }$ at distance $z$ from $O^{\prime }$ such that 
\begin{equation*}
|O^{\prime }P^{\prime }|\text{ }=|O^{\prime }N^{\prime }|\text{ }+|N^{\prime
}P^{\prime }|\text{ }=\delta z+\delta (90-p).
\end{equation*}%
Since the distance of the pole $N$ to the equator is $90^{0},$ the point $%
N^{\prime }$ can be considered as the image of the pole $N$ by the
projection $f.$

The distance $|O^{\prime }N^{\prime }|$ measures
how far is the point $O^{\prime }$ beyond the pole $N^{\prime }.$
Substituting in the previous equation $|O^{\prime }P^{\prime }|$ from (\ref%
{1}) we get%
\begin{equation}
z=\frac{(q-p)\cos p}{\cos p-\cos q}-90+p.  \tag{4}  \label{4}
\end{equation}

If $A^{\prime }=f(A)$ and $B^{\prime }=f(B)$ we have  $|A^{\prime }O^{\prime }|$ $=\delta (90^{\mathrm{o}}-a+z).$ Multiplying this quantity
by $\omega $ and using (\ref{3}) we obtain that the length of the circular
arc of radius $\vert O'A'\vert$ ad angle $\omega$ is approximately $\vert A'A'_A\vert$, and, using \ref{3}, $A'A'_A\vert$ is 

\begin{equation*}
\frac{\delta (90^{\mathrm{o}}-a+z)(\cos p-\cos q)}{\alpha (q-p)}.
\end{equation*}%
Since the length of one degree on the parallel passing through $A$ in $S$ is 
$\delta \cos a$ the difference 
\begin{equation*}
\frac{\delta (90^{\mathrm{o}}-b+z)(\cos p-\cos q)}{\alpha (q-p)}-\delta \cos a
\end{equation*}%
shows the error of the projection at the extremity $A.$

 Similarly the error of the projection at the extremity $B$ is 
\begin{equation*}
\frac{\delta (90^{\mathrm{o}}-p+z)(\cos p-\cos q)}{\alpha (q-p)}-\delta \cos b.
\end{equation*}
  Since the points $P$ and $Q$ belong to the meridian $AB$ and since they lie between $A$ and $B$, and since
 the errors at the two extremities $A$ and $B$ are assumed to be equal to
each other, we obtain the equation 
$$
\frac{\delta (90^{\mathrm{o}}-a+z)(\cos p-\cos q)}{\alpha (q-p)}-\delta \cos a$$
$$=\frac{%
\delta (90^{\mathrm{o}}-b+z)(\cos p-\cos q)}{\alpha (q-p)}-\delta \cos b
$$
which can take the form%
\begin{equation}
(a-b)(\cos p-\cos q)=(q-p)(\cos a-\cos b).  \tag{5}  \label{5}
\end{equation}%
Considering $a=40^{\mathrm{o}},$ $b=70^{\mathrm{o}},$ $p=50^{\mathrm{o}},$ $q=60^{\mathrm{o}}$ we may verify
that the two members of (\ref{5}) are almost equal. But in general, given $a$
and $b$, it is difficult to detect values $p$ and $q$ inside $(a,b)$ such that (%
\ref{5}) is satisfied.

Thus, instead of the quantities $p$ and $q$ it is easier to work with the
quantity $z$ which expresses the displacement of the point $O^{\prime }$ from $%
N^{\prime }.$ Below we shall express $z$ as a function of $a$ and $b$ and so $z$ can be calculated without using \ref{4}. One degree of parallel at the extremity $A^{\prime }$ will be
equal to $\alpha \delta (90-a+z)\omega $ while the measure of one degree of the parallel of $S$ passing trough $A$ is $\delta \cos a.$ Therefore the
error at $A$ is%
\begin{equation*}
\alpha \delta (90-a+z)\omega -\delta \cos a.
\end{equation*}%
Similarly at the extremity $B$ the error is 
\begin{equation*}
\alpha \delta (90-b+z)\omega -\delta \cos b.
\end{equation*}%
Setting that these errors are equal among themselves we obtain the equation 
\begin{equation*}
\alpha (a-b)\omega =\cos a-\cos b
\end{equation*}%
from which we get%
\begin{equation}
\omega =\frac{\cos a-\cos b}{\alpha (a-b)}. \tag{6}  \label{6}
\end{equation}%

Now, after making equal the errors of the projection at the extremities $A$
and $B$ we impose the following additional condition:

\begin{itemize}
\item The error at $A$ and $B$ is equal to the one that occurs at
the midpoint $X$ of the interval $AB$ whose latitude is $\frac{a+b}{2}$, which is supposed to be the maximal error of the projection.
\end{itemize}

The error at $X$ is 
\begin{equation*}
\alpha \delta (90-\frac{a+b}{2}+z)\omega -\delta \cos \frac{a+b}{2}.
\end{equation*}%
Note that we have to assume that the sign of this error is negative so it is
necessary to put the error equal to 
\begin{equation*}
\delta \cos \frac{a+b}{2}-\alpha \delta (90-\frac{a+b}{2}+z)\omega .
\end{equation*}%
Since we have assumed that the errors at $A,$ $B$ and $X$ are equal we have
the following equations

\begin{equation*}
\alpha \delta (90-a+z)\omega -\delta \cos a=\delta \cos \frac{a+b}{2}-\alpha
\delta (90-\frac{a+b}{2}+z)\omega
\end{equation*}%
and%
\begin{equation*}
\alpha \delta (90-b+z)\omega -\delta \cos b=\delta \cos \frac{a+b}{2}-\alpha
\delta (90-\frac{a+b}{2}+z)\omega .
\end{equation*}%
Substituting now the value of $\omega $ from (\ref{6}) in one or the other
of the two previous equations we find the equation 
\begin{equation}
\frac{\alpha (180-\frac{3}{2}a-\frac{1}{2}b+2z)(\cos a-\cos b)}{b-a}=\cos
a+\cos \frac{a+b}{2}.  \tag{7}  \label{7}
\end{equation}%
From (\ref{7}) we can determine $z$ as a function of $a$ and $b.$

\section{Applications of Delisle's projection to the map of the Russian
Empire}

In this section,\index{map!Delisle}\index{projection!Delisle}\index{Delisle map}\index{map!Russian Empire}\index{general map of the Russian Empire} we will apply the previous analysis to the case of the map
of the Russian Empire. Taking $a=40^{\mathrm{o}}$ and $b=70^{\mathrm{o}}$ we have $\frac{%
a+b}{2}=55^{\mathrm{o}}$ and therefore from (\ref{6}) we obtain  
\begin{equation}
\omega =\frac{\cos 40^{\mathrm{o}}-\cos 70^{\mathrm{o}}}{30\alpha }=\frac{0,4240243}{0,5235987}
\tag{8}  \label{8}
\end{equation}%
in radians and hence (in degrees) $\omega =48^{\prime }44^{\prime \prime }.$ Thus, from (%
\ref{7}) we have%
\begin{equation*}
\alpha (85^{\mathrm{o}}+2z)\omega =\cos 40^{\mathrm{o}}+\cos 55^{\mathrm{o}}=1,33962.
\end{equation*}%
Then, substituting the values of $\alpha $ and $\omega $ we get 
\begin{equation*}
85^{\mathrm{o}}+2z=\frac{1,33962}{0,0141}=95^{\mathrm{o}}\Rightarrow z=5^{\mathrm{o}}.
\end{equation*}

In the previous section we have assumed that the maximal error occurs at
the midpoint of the arc $AB.$ However, the point of maximal error could
deviate from the midpoint. So, we will find exactly the point $X$ at which
the maximal error occurs. For this we consider the function 
\begin{equation*}
e:[a,b]\rightarrow \mathbb{R}
\end{equation*}%
with 
\begin{equation}
e(x)=\omega \alpha  (90^{\mathrm{o}}-x+z)-\cos x,  \tag{9}  \label{9}
\end{equation}%
where $x$ is expressed in degrees.
Then we have $e^{\prime }(x)=\alpha\sin x-\alpha \omega\alpha ,$ $e^{\prime \prime
}(x)=\alpha^2\cos x>0$ in $[a,b].$ Therefore the maximum occurs for $x$ such that $%
\sin x=\omega ,$ where $\omega $ is given by (\ref{8}). Therefore we have
\begin{equation*}
\sin x=\frac{0,4240243}{0,5235987}
\end{equation*}%
and we get $x=54^{\mathrm{o}}4^{\prime },$ which is a point that differs little
from the midpoint of the arc $AB.$

Having found the above value for $x,$ the error at this point will be 
\begin{equation*}
\alpha (90^{\mathrm{o}}-x+z)\omega -\cos x
\end{equation*}%
and if we impose that this quantity, taken with a negative sign, is equal to
the error at $A$ we have
\begin{equation*}
\cos x-\alpha (90^{\mathrm{o}}-x+z)\omega =\alpha (90^{\mathrm{o}}-a+z)\omega -\cos a
\end{equation*}%
and thus we get the equation%
\begin{equation*}
\alpha (180^{\mathrm{o}}-a-x+2z)\omega =\cos x+\cos a
\end{equation*}%
from which the value $z$ can be again drawn out. More precisely, since $x=54^{\mathrm{o}}4'$ the equation takes the form%
\begin{equation*}
85\frac{14}{15}+2z=\frac{\cos a+\cos x}{\alpha \omega }=95^{\mathrm{o}}56^{\prime };
\end{equation*}%
therefore $z=5^{\mathrm{o}}$ since $\omega =0,8098270$ (which corresponds in degrees
to $\omega =48^{\mathrm{o}}44^{\prime }).$

Replacing the value $z=5^{\mathrm{o}}$ in (\ref{9}), the function $e(x)=\alpha \omega
(90^{\mathrm{o}}-x+z)-\cos x$ takes the form $e(x)=\alpha \omega (95^{\mathrm{o}}-x)-\cos x.$
This last function is almost zero for $x=50^{0}$ or $x=60^{0}.$ Therefore,
if we consider the parallels defined for latitudes $p=50^{0}$ and $q=60^{0}$
we deduce that along these parallels the ratio of the degrees of latitude to
the degrees of longitude is almost constant. On the other hand, finding the
values of $x$ which are roots of the equation $\alpha \omega
(95^{\mathrm{o}}-x)-\cos x=0,$ we may find exactly the latitudes of parallels along
which the ratio of the degrees of latitude to the degrees of longitude is
constant. In this way one can answer Euler's second question mentioned in
the introduction.\smallskip

\noindent \textbf{Note}: We have seen that the function $e(x)=\alpha \omega
(90^{\mathrm{o}}-x+z)-\cos x$ has a unique critical point, in particular a minimum at
a point, say at $x_{0}.$ Since $e(a)=e(b)=|e(x_{0})|$ it follows that $%
e(a)=e(b)=-e(x_{0})$ and that the equation $\alpha \omega (95^{\mathrm{o}}-x)-\cos
x=0 $ has exactly two roots.\smallskip

Now let us compute how big is the error at $A$ and $B.$ The error at $A$ is%
\begin{equation*}
\alpha \omega (90^{\mathrm{o}}-a+z)-\cos a=55\alpha \omega -0,7660444.
\end{equation*}%
Since $\alpha \omega =0,01410$ the error turns out to be equal to $0,00946$
and since this error is expressed in fractions of a meridian degree,
assigning to such a degree the length of $15$ miles the measure of the error
is $0,14190$ miles, that is the seventh part of a mile. Therefore at the
extremity $B,$ where the latitude is $70^{\mathrm{o}},$ and hence one degree of
parallel is equal to $0,34202,$ the error is equal to the thirty-eight part
of a mile, which is easily tolerated.

Finally, it is quite easy now to construct Delisle's\index{map!Delisle}\index{projection!Delisle}\index{Delisle map}\index{Delisle projection} map of the Russian
Empire, without considering the points $P^{\prime }$ and $Q^{\prime }$ of
Paragraph 2. For this, we consider first a segment $A^{\prime }B^{\prime }$
in $E^{2}$ which will be the image of the arc $AB$ of the principal meridian. At a
distance of $5$ degrees from $B^{\prime }$ (or $5\delta )$ we consider the
point $O^{\prime }$ on the extension of the segment $A^{\prime }B^{\prime }.$
In the following, at the point $A^{\prime }$ which is at a distant equal to $%
55^{\mathrm{o}}$ from $O^{\prime },$ we consider with center $O^{\prime }$ the circle 
$S_{A^{\prime }}$ of radius $O^{\prime }A^{\prime }$ and this circle will be
the image of the parallel passing through $A.$ On $S_{A^{\prime }}$ a degree
of longitude will be equal to $55\alpha \omega =0.77550,$ so the division on
this circle can be completed and all the meridians are easily drawn as
rays emanating from $O^{\prime }$ and passing from the points
of the subdivision of $S_{A^{\prime }}.$ The images of the parallels of $S$ whose
latitude differs by one degree will be circles of center $O^{\prime }$ which
have distance $\delta $ among them.

\end{document}